\documentclass{amsart}
\usepackage{amssymb,latexsym}

\theoremstyle{plain}
\newtheorem{theorem}{Theorem}

\theoremstyle{definition}
\newtheorem{definition}{Definition}

\theoremstyle{remark}

\numberwithin{equation}{section}

\begin{document}
\title[Some problems in differentiation]
{Some problems in differentiation}
\author{Warren P.~Johnson}
\address{Connecticut College\\
270 Mohegan Avenue\\
New London, CT 06320}
\email{wpjoh@conncoll.edu}
\keywords{derivatives, inverse functions, implicit functions, parametric
equations}

\maketitle

\section{Introduction}\label{S:intro}

In this note we discuss three problems.

\textbf{Problem 1:}  What is the $n^{\text{th}}$ derivative of an
inverse function?

\textbf{Problem 2:}  What is the $n^{\text{th}}$ derivative of a
function given parametrically?

\textbf{Problem 3:}  What is the $n^{\text{th}}$ derivative of an
implicit function?

Let us be more precise.  Suppose $x=f(t)$ and $y=g(t)$ for nice 
enough functions $f$ and $g$.  What is $d^ny/dx^n$ in terms of
the derivatives of $f$ and $g$?  Or suppose that an equation
$F(x,y)=0$ defines $y$ as an implicit function of $x$, and that
$F$ is nice enough to have equality of all the relevant mixed
partial derivatives.  What is $d^ny/dx^n$ in terms of the partial 
derivatives of $F$?  These are problems 2 and 3 respectively.

Note that problem 1 is the special case $g(t)=t$ of problem 2,
and the special case $F(x,y)=f(y)-x$ of problem 3, so
it does not need a separate treatment.  Moreover, I have discussed 
it before in \cite{johnson}, which contains my rediscovery of
Sylvester's solution \cite{sylvester}.  Problem 1 was apparently
solved first by Murphy \cite{murphy}.

Sylvester pointed out in a postscript to \cite{sylvester} that his 
argument extends without great difficulty to problem 2, but he
did not give the details.  A survey of the problem was given later 
by Gallop \cite{gallop}, who had read \cite{sylvester} but was 
not sympathetic to it.  The following treatment is, I believe,
close to what Sylvester had in mind.

We begin by working out the first several derivatives.  Evidently we have 
\[
\frac{dy}{dx}=\frac{\frac{dy}{dt}}{\frac{dx}{dt}}=\frac{g'(t)}{f'(t)}=
g'(t)\left[f'(t)\right]^{-1}.
\]
Then
\begin{align*}
\frac{d^2y}{dx^2}&=\frac{\frac{d}{dt}\left(\frac{dy}{dx}\right)}{\frac{dx}{dt}}
=\frac1{f'(t)}\left\{g''(t)\left[f'(t)\right]^{-1}-
g'(t)\left[f'(t)\right]^{-2}f''(t)\right\}\\
&=g''(t)\left[f'(t)\right]^{-2}-g'(t)f''(t)\left[f'(t)\right]^{-3}
\end{align*}
and 
\begin{equation}\label{E:three}
\begin{aligned}
\frac{d^3y}{dx^3}&=g'''(t)\left[f'(t)\right]^{-3}-3g''(t)f''(t)\left[f'(t)
\right]^{-4}\\
&\quad-g'(t)f'''(t)\left[f'(t)\right]^{-4}+3g'(t)\left[f''(t)\right]^2
\left[f'(t)\right]^{-5}.
\end{aligned}
\end{equation}

To illustrate the combinatorics behind these formulas we take
(\ref{E:three}).  Besides the order $n$ of the derivative (here
$n=3$) there is another crucial parameter that we will call $k$.
Each term of (\ref{E:three}) represents a partition of 
$\{1,2,\dots,k+3\}$ into $k+1$ blocks, where $k+3$ is the exponent
of $\left[f'(t)\right]^{-1}$ and the sign of each term is $(-1)^k$.
The term $g'''(t)\left[f'(t)\right]^{-3}$, with $k=0$, represents
$\{1,2,3\}$.  The term $-3g''(t)f''(t)\left[f'(t)\right]^{-4}$,
with $k=1$, represents the three partitions $\{1,2\},\{3,4\}$ and
$\{1,3\},\{2,4\}$ and $\{1,4\},\{2,3\}$, where the block containing
$1$ and one other element corresponds to $g''(t)$ and the other
block to $f''(t)$.  The block containing $1$ always corresponds
to and has the same size as the order of the $g$ derivative, so
the other $k=1$ term $-g'(t)f'''(t)\left[f'(t)\right]^{-4}$ 
corresponds to $\{1\},\{2,3,4\}$.  The remaining term 
$3g'(t)\left[f''(t)\right]^2\left[f'(t)\right]^{-5}$ has $k=2$
and represents the partitions $\{1\},\{2,3\},\{4,5\}$ and
$\{1\},\{2,4\},\{3,5\}$ and $\{1\},\{2,5\},\{3,4\}$.

Note that there is never an $f'(t)$ term in the numerators of these
formulas; only higher derivatives of $f$ can occur.  This means that
the corresponding partitions cannot have any singleton blocks besides
$\{1\}$, which in turn explains why $k$ can't exceed 2 in this example:
if $k=3$, we would have to put $\{1,2,3,4,5,6\}$ into 4 blocks without
putting any of $\{2,3,4,5,6\}$ in a singleton, which is impossible.

Let us try to predict $d^4y/dx^4$ from this point of view.  At the
same time we illustrate a definition we shall need for the general
case.  We have to consider partitions of $\{1,2,\dots,k+4\}$ with 
$k+1$ blocks, where only 1 can be in a singleton.  If $k=0$ this 
means four elements and one block.  The only such partition is
$\{1,2,3,4\}$, which gives us the term $g''''(t)\left[f'(t)\right]^{-4}$.
If we set $\boldsymbol{P}_{4,0}(t)=g''''(t)$ then this term is
$\left[f'(t)\right]^{-4}\boldsymbol{P}_{4,0}(t)$.

When $k=1$ we have five elements and two blocks, so the block sizes
are either one and four or two and three.  With one and four we can
only have $\{1\},\{2,3,4,5\}$, which corresponds to $-g'(t)f''''(t)
\left[f'(t)\right]^{-5}$.  With two and three there are two possibilities:
1 is in a block of size two, or 1 is in a block of size three.  In the 
first case there are 4 other elements that could be together with 1,
so these partitions give the term $-4g''(t)f'''(t)\left[f'(t)\right]^{-5}$.
In the second case there are $\binom 42=6$ ways to choose the block
containing 1, so these partitions give the term $-6g'''(t)f''(t)
\left[f'(t)\right]^{-5}$.  If we set $\boldsymbol{P}_{4,1}(t)=
g'(t)f''''(t)+4g''(t)f'''(t)+6g'''(t)f''(t)$, then these three terms
together are $-\left[f'(t)\right]^{-5}\boldsymbol{P}_{4,1}(t)$.

If $k=2$ we have six elements and three blocks.  Since there can only
be one singleton, the block sizes must be either 3-2-1 or 2-2-2.  In
the former case we have the block $\{1\}$, and we complete the partition
by choosing the block of size two in any of $\binom 52=10$ ways.
Therefore these partitions correspond to $10g'(t)f''(t)f'''(t)
\left[f'(t)\right]^{-6}$.  If all three blocks have size two, then
we pick an element to put with 1, and then pick one to put with the
smallest remaining element, and the partition is determined, so there
are $5\cdot3=15$ partitions in this case and the corresponding term is
$15g''(t)\left[f''(t)\right]^2\left[f'(t)\right]^{-6}$.  If we set 
$\boldsymbol{P}_{4,2}(t)=10g'(t)f''(t)f'''(t)+15g''(t)\left[f''(t)
\right]^2$ then these two terms together are $\left[f'(t)\right]^{-6}
\boldsymbol{P}_{4,2}(t)$.

If $k=3$ we have seven elements and four blocks.  This forces at
least one singleton, but there can't be more than one, so the
only possibility is to have the block $\{1\}$ and the other six
elements in three blocks of size two.  As above there are 15 such
partitions, so this gives the term $-15g'(t)\left[f''(t)\right]^3
\left[f'(t)\right]^{-7}$; or, setting $\boldsymbol{P}_{4,3}(t)=
15g'(t)\left[f''(t)\right]^3$, the term $-\left[f'(t)\right]^{-7}
\boldsymbol{P}_{4,3}(t)$. It is also clear from this example that
we would get an impossible number of singleton blocks if $k>3$,
so we define $\boldsymbol{P}_{4,k}(t)=0$ in this case. 

Adding all these terms together we have
\[
\frac{d^4y}{dx^4}=\sum_{k=0}^3(-1)^k\left[f'(t)\right]^{-4-k}
\boldsymbol{P}_{4,k}(t),
\]
where we could just as well leave the upper limit of the sum unrestricted.
When written out broadly this says
\begin{align*}
\frac{d^4y}{dx^4}&=g''''(t)\left[f'(t)\right]^{-4}-6g'''(t)f''(t)\left[f'(t)\right]^{-5}
+15g''(t)\left[f''(t)\right]^2\left[f'(t)\right]^{-6}\\
&\quad-4g''(t)f'''(t)\left[f'(t)\right]^{-5}
+10g'(t)f''(t)f'''(t)\left[f'(t)\right]^{-6}\\
&\quad-g'(t)f''''(t)\left[f'(t)\right]^{-5}-15g'(t)\left[f''(t)\right]^3
\left[f'(t)\right]^{-7}
\end{align*}
which is right.

To state the general result we make the definition illustrated above.

\begin{definition}\label{D:partsum} Given two functions $f(t)$ and
$g(t)$ with at least $n$ derivatives, define $\boldsymbol{P}_{n,k}(t)$ 
as the sum over all the partitions of $\{1,2,\dots,k+n\}$ with $k+1$ 
blocks where only the block containing 1 can be a singleton, and as
above each partition corresponds to a product of derivatives, where 
the size of the block containing 1 is the order of the derivative of 
$g$ and the sizes of the other blocks are orders of derivatives of $f$.
\end{definition}

We note some extreme cases.  If $k=0$ we have only the partition
$\{1,2,\dots,n\}$, so $\boldsymbol{P}_{n,0}(t)=g^{(n)}(t)$ for $n\ge1$.
At the other end, if $k=n-1$ we have partitions of $\{1,2,\dots,2n-1\}$
into $n$ blocks, where only the element 1 can be in a singleton.
Therefore all the other blocks must be doubletons, and there are
$1\cdot3\cdot5\cdots(2n-3)$ ways to choose them, so
\[
\boldsymbol{P}_{n,n-1}(t)=1\cdot3\cdot5\cdots(2n-3)g'(t)\left[f''(t)
\right]^{n-1} \quad\text{for $n\ge1$.}
\]
If $k>n-1$ then we have an impossible number of singletons, so
$\boldsymbol{P}_{n,k}(t)=0$ if $k>n-1$, or for that matter if $k<0$.
It is also not hard to see that
\begin{multline*}
\boldsymbol{P}_{n,n-2}(t)=1\cdot3\cdot5\cdots(2n-3)g''(t)\left[f''(t)
\right]^{n-2}\\
+\binom{2n-3}{3}1\cdot3\cdot5\cdots(2n-7)g'(t)\left[f''(t)
\right]^{n-3}f'''(t) \quad\text{for $n\ge2$.}
\end{multline*}

We need the following recurrence for $\boldsymbol{P}_{n,k}(t)$.

\begin{equation}\label{E:recur}
\boldsymbol{P}_{n+1,k}(t)=\boldsymbol{P}'_{n,k}(t)+(n+k-1)f''(t)
\boldsymbol{P}_{n,k-1}(t).
\end{equation}

We can obtain partitions of $\{1,2,\dots,k+n+1\}$ with $k+1$ blocks
from a partition of $\{1,2,\dots,k+n\}$ with $k+1$ blocks by adding
the element $k+n+1$ to each existing block in turn.  This corresponds
perfectly to how $d/dt$ acts on each term of the sum, increasing the
order of one derivative by one while leaving the others alone, and
doing this for each factor in turn.  The only problem is that we
cannot obtain \emph{all} of the partitions of the desired type this
way:  we are missing the ones where $k+n+1$ is in a doubleton with
one of $\{2,3,\dots,k+n\}$, because those elements could not have
been in singleton blocks.  These partitions come from the term
$(n+k-1)f''(t)\boldsymbol{P}_{n,k-1}(t)$, for we know that 
$\boldsymbol{P}_{n,k-1}(t)$ corresponds to partitions of 
$\{1,2,\dots,k+n-1\}$ with $k$ blocks.  We can add
$\{k+n,k+n+1\}$ as a doubleton to each of these, and then we can
switch $k+n$ with any of $\{2,3,\dots,k+n-1\}$ to get the remaining
ones.  This proves (\ref{E:recur}), and now we can prove the following
formula.

\begin{theorem}\label{T:parametric} If $x=f(t)$ and $y=g(t)$, where
$f$ and $g$ have at least $n$ derivatives, then
\begin{equation}\label{E:parametric}
\frac{d^ny}{dx^n}=\sum_{k=0}^{n-1}(-1)^k\left[f'(t)\right]^{-n-k}
\boldsymbol{P}_{n,k}(t).
\end{equation}
\end{theorem}

We have checked the cases $n=1,2,3,4$ of this.  Suppose it holds
for $n$, and take $d/dx$ of the right side of (\ref{E:parametric}) by
taking $d/dt$ and dividing by $dx/dt=f'(t)$. 
It is convenient to leave the sum unrestricted, as we may since
$\boldsymbol{P}_{n,k}(t)=0$ if $k<0$ or $k>n-1$.  We get
\begin{multline*}
\frac1{f'(t)}\left[\sum_j(-1)^j(-n-j)\left[f'(t)\right]^{-n-j-1}f''(t)
\boldsymbol{P}_{n,j}(t)+\sum_j(-1)^j\left[f'(t)\right]^{-n-j}
\boldsymbol{P}'_{n,j}(t)\right]\\
=\sum_j(-1)^{j+1}(n+j)\left[f'(t)\right]^{-n-1-(j+1)}f''(t)
\boldsymbol{P}_{n,j}(t)+\sum_j(-1)^j\left[f'(t)\right]^{-n-1-j}
\boldsymbol{P}'_{n,j}(t)\\
=\sum_k(-1)^k(n+k-1)\left[f'(t)\right]^{-n-1-k}f''(t)
\boldsymbol{P}_{n,k-1}(t)+\sum_k(-1)^k\left[f'(t)\right]^{-n-1-k}
\boldsymbol{P}'_{n,k}(t).
\end{multline*}

Because of (\ref{E:recur}) we know that this is
\[
\sum_k(-1)^k\left[f'(t)\right]^{-n-1-k}\boldsymbol{P}_{n+1,k}(t).
\]
Therefore (\ref{E:parametric}) holds for $n+1$ if it holds for $n$,
and hence Theorem \ref{T:parametric} is true.

Next we take up problem 3.  It has been treated in two papers by Comtet
\cite{comtet1}, \cite{comfio}, whose excellent book \cite{comtet2} also has
references (p.~153) to some older literature.  More recently it has been
discussed in \cite{wilde} and \cite{nahay}.  The forthcoming paper by Shaul
Zemel \cite{zemel}
gives a valuable survey of the various approaches, including my own below,
and a new method.   I am putting this version of my paper on the arXiv at
his request.

If $F(x,y)=0$ then we have
\[
\frac{{\partial}F}{{\partial}x}+\frac{{\partial}F}{{\partial}y}\,
\frac{dy}{dx}=0,
\]
so
\[
\frac{dy}{dx}=-\frac{\frac{{\partial}F}{{\partial}x}}{\frac{{\partial}F}{{\partial}y}},
\]
or, as we will prefer to write,
\[
\frac{dy}{dx}=-F_xF_y^{-1}.
\]
Differentiating this with respect to $x$ we have
\begin{align*}
\frac{d^2y}{dx^2}&=-\frac{\partial}{{\partial}x}\,F_xF_y^{-1}
-\frac{dy}{dx}\,\frac{\partial}{{\partial}y}\,F_xF_y^{-1}\\
&=-F_{xx}F_y^{-1}+F_xF_y^{-2}F_{xy}+F_xF_y^{-1}\left(F_{xy}F_y^{-1}
-F_xF_y^{-2}F_{yy}\right)\\
&=-F_{xx}F_y^{-1}+2F_xF_{xy}F_y^{-2}-F_x^2F_{yy}F_y^{-3}.
\end{align*}
As before, we want to associate a family of set partitions to these
formulas.  The general result we are after has the form
\[
\frac{d^ny}{dx^n}=\sum_{k\ge1}(-1)^kF_y^{-k}\boldsymbol{I}_{n,k},
\]
where $\boldsymbol{I}_{n,k}$ has still to be explained.  We have
the examples $\boldsymbol{I}_{2,1}=F_{xx}$, $\boldsymbol{I}_{2,2}
=2F_xF_{xy}$, $\boldsymbol{I}_{2,3}=F_x^2F_{yy}$, and 
$\boldsymbol{I}_{2,k}=0$ for $k>3$.  In general $\boldsymbol{I}_{n,k}$
will be a sum over partitions of $\{1,2,\dots,n+k-1\}$ with $k$ blocks.
Let us call $\{1,2,\dots,n\}$ the \emph{small} elements, and
$\{n+1,\dots,n+k-1\}$ the \emph{large} elements.  The small elements
will correspond to $x$ derivatives and the large elements to $y$
derivatives, and the partitions are restricted only in that the large
elements can't be in singleton blocks.  When $k=1$ all the elements
are small and must be together, giving us the term $F_{xx}$ when $n=2$
and in general the term ${\partial}^nF/{\partial}x^n$.  When $n=2=k$ 
we have partitions of $\{1,2,3\}$ into two blocks, where 3 is large
and hence can't be in a singleton.  There are two
such partitions, $\{1\},\{2,3\}$ and $\{1,3\},\{2\}$, and with
1 and 2 corresponding to $x$ and 3 to $y$ these give the term
$2F_xF_{xy}$.  When $n=2$ and $k=3$ we have partitions of $\{1,2,3,4\}$
with three blocks, where 3 and 4 are large.  There is
only one such partition, $\{1\},\{2\},\{3,4\}$, and it corresponds
to $F_x^2F_{yy}$.  If $n=2$ and $k=4$ we would have to put $\{1,2,3,4,5\}$
into 4 blocks without putting 3, 4, or 5 into a singleton, which is
impossible, and similarly there are no more partitions with $n=2$
and a larger $k$.

Let us work out $d^3y/dx^3$ from this point of view;
we have $n=3$ and various values of $k$.
As before, when $k=1$ we have only the term 
$\boldsymbol{I}_{3,1}=F_{xxx}$.  When
$k=2$ we have partitions of $\{1,2,3,4\}$ with two blocks, where 
4 can't be in a singleton.  All three partitions with two blocks
of size two are admissible, and they give us the term $3F_{xx}F_{xy}$.
Of the four partitions with block sizes one and three we can only
use those with 4 in a triplet, and they give us the term $3F_xF_{xxy}$.  
Therefore $\boldsymbol{I}_{3,2}=3F_{xx}F_{xy}+3F_xF_{xxy}$.

When $n=3=k$ we have partitions of $\{1,2,3,4,5\}$ with three blocks,
where 4 and 5 are large.  The block sizes must be either 3-1-1 or
2-2-1.  In the former case the singletons must be two of $\{1,2,3\}$,
so there are only three possibilities, which together give the term
$3F_x^2F_{xyy}$.  In the latter case there are two subcases.  If
the two large elements 4 and 5 are together in a doubleton, the
partition is determined by which of the three small elements is
alone, and we get the term $3F_xF_{xx}F_{yy}$.  If 4 and 5 are
in different doubletons then there are six possibilities, and
we get the term $6F_xF_{xy}^2$.  Therefore $\boldsymbol{I}_{3,3}=
3F_x^2F_{xyy}+3F_xF_{xx}F_{yy}+6F_xF_{xy}^2$.

When $n=3$ and $k=4$ we have partitions of $\{1,2,3,4,5,6\}$ with
four blocks, where 4,5,6 are large.  The block sizes must be
either 3-1-1-1 or 2-2-1-1.  There is only one admissible partition
in the former case, namely $\{1\},\{2\},\{3\},\{4,5,6\}$, which
corresponds to the term $F_x^3F_{yyy}$.  In the latter case we
can choose one of 1,2,3 to be in a doubleton with one of 4,5,6,
and after that the partition is determined (since the other two
small elements must be the singletons), so there are nine possibilities, 
which together contribute the term $9F_x^2F_{xy}F_{yy}$.
Hence $\boldsymbol{I}_{3,4}=F_x^3F_{yyy}+9F_x^2F_{xy}F_{yy}$.

When $n=3$ and $k=5$ we have partitions of $\{1,2,3,4,5,6,7\}$
with five blocks, where 4,5,6,7 are large.  With this constraint
we cannot have block sizes 3-1-1-1-1, so we have only to consider
the case 2-2-1-1-1.  The singletons must be 1,2,3, and there are
three such partitions, which collectively give us $3F_x^3F_{yy}^2=
\boldsymbol{I}_{3,5}$.

In general there are no partitions of the desired type when $k>2n-1$,
so we have all of them for $n=3$.  If $k=2n$ then we would have to
put $3n-1$ elements into $2n$ blocks, which forces at least $n+1$
singletons; but only the $n$ small elements can be in singletons,
so this is impossible.  If $k>2n$ then even more elements have to 
be in singletons, which is even more impossible.  Therefore we define
$\boldsymbol{I}_{n,k}=0$ if $k>2n-1$.

The recurrence we will need for $\boldsymbol{I}_{n,k}$ is rather
complicated.

\begin{equation}\label{E:recurr}
\boldsymbol{I}_{n+1,k}=\frac{\partial}{{\partial}x}\,\boldsymbol{I}_{n,k}
+F_x\,\frac{\partial}{{\partial}y}\,\boldsymbol{I}_{n,k-1}+(k-1)F_{xy}\,\boldsymbol{I}_{n,k-1}+(k-2)F_xF_{yy}\,\boldsymbol{I}_{n,k-2}.
\end{equation}

The left side is a sum over partitions of $\{1,2,\dots,n+k\}$ with $k$
blocks, where $\{1,2,\dots,n+1\}$ are small and $\{n+2,\dots,n+k\}$ are
large, and we have to argue that the right side also represents this
class of partitions.  As before, a derivative adds an element to an
existing block while leaving the others alone---here an $x$
derivative adds a small element and a $y$ derivative a large one.

We interpret the first term $\frac{\partial}{{\partial}x}\,
\boldsymbol{I}_{n,k}$ as adding the element 1 to an existing
block in a partition of $\{1,2,\dots,n+k-1\}$ with $k$ blocks
of the desired type and relabeling all the other elements up
one; thus we have one more small element and the same number
of large elements, as desired.  Note that with this operation
1 can't be in a singleton block, nor can it be in a doubleton
with a large element, but it gives all the partitions of the
desired type except for these restrictions.

The second term $F_x\,\frac{\partial}{{\partial}y}\,
\boldsymbol{I}_{n,k-1}$ should be an operation on partitions
of $\{1,2,\dots,n+k-2\}$ with $k-1$ blocks where $\{1,2,\dots,n\}$
are small.  We think of it as adding $\{1\}$ as a singleton block,
relabeling all the other elements up one, and adding $n+k$ to an
existing block.  This gives the correct numbers of small and large
elements, but in addition to putting 1 in a singleton it also
means that $n+k$ is not in a doubleton with another large element.

The third term $(k-1)F_{xy}\,\boldsymbol{I}_{n,k-1}$ is an 
operation on the same partitions as the second term.  We
think of it as first adding $\{1,n+k\}$ as a doubleton and
relabeling all the other elements up one, then switching
$n+k$ with each of the $k-2$ relabeled large elements 
$\{n+2,\dots,n+k-1\}$ in turn to create $k-2$ more partitions
of the desired type.  These are the missing ones from the
first term, having 1 in a doubleton with a large element.

Finally, the last term $(k-2)F_xF_{yy}\,\boldsymbol{I}_{n,k-2}$
should be an operation on partitions of $\{1,2,\dots,n+k-3\}$
of the desired type with $k-2$ blocks where $\{1,2,\dots,n\}$ 
are small.  Here we first add $\{1\}$ as a singleton block, 
relabel the other elements up one, and add the doubleton $\{n+k-1,n+k\}$.
Then we switch $n+k-1$ with each of the $k-3$ relabeled large elements
$\{n+2,\dots,n+k-2\}$ in turn to create $k-3$ more partitions of
the desired type.  These are the missing ones from the second term,
where 1 is in a singleton and $n+k$ is in a doubleton with another
large element.  This proves (\ref{E:recurr}).

\begin{theorem}\label{T:implicit}
If $F(x,y)=0$ for a sufficiently nice function $F$, then
\begin{equation}\label{E:implicit}
\frac{d^ny}{dx^n}=\sum_{k=1}^{2n-1}(-1)^kF_y^{-k}\boldsymbol{I}_{n,k}
\end{equation}
with $\boldsymbol{I}_{n,k}$ as defined above.
\end{theorem}

In the proof it is convenient to leave the sum unrestricted, as we
may since $\boldsymbol{I}_{n,k}=0$ if $k<1$ or if $k>2n-1$.  On the
right side of (\ref{E:implicit}) we calculate the derivative with 
respect to $x$ as
\[
\frac{\partial}{{\partial}x}+\frac{dy}{dx}\,\frac{\partial}{{\partial}y}=
\frac{\partial}{{\partial}x}-\frac{F_x}{F_y}\,\frac{\partial}{{\partial}y}.
\]

Applying $\frac{\partial}{{\partial}x}$ to the right side of 
(\ref{E:implicit}) we get

\begin{multline*}
\sum_k(-1)^k(-k)F_y^{-k-1}F_{xy}\boldsymbol{I}_{n,k}+\sum_k(-1)^kF_y^{-k}\,\frac{\partial}{{\partial}x}\,\boldsymbol{I}_{n,k}\\
=\sum_j(-1)^{j+1}jF_y^{-(j+1)}F_{xy}\boldsymbol{I}_{n,j}+
\sum_j(-1)^jF_y^{-j}\,\frac{\partial}{{\partial}x}\,\boldsymbol{I}_{n,j}\\
=\sum_k(-1)^k(k-1)F_y^{-k}F_{xy}\boldsymbol{I}_{n,k-1}+\sum_k(-1)^kF_y^{-k}\,\frac{\partial}{{\partial}x}\,\boldsymbol{I}_{n,k}.
\end{multline*}

Applying $-\frac{F_x}{F_y}\,\frac{\partial}{{\partial}y}$ to the right 
side of (\ref{E:implicit}) we get

\begin{multline*}
-\frac{F_x}{F_y}\,\sum_k(-1)^k(-k)F_y^{-k-1}F_{yy}\boldsymbol{I}_{n,k}
-\frac{F_x}{F_y}\,\sum_k(-1)^kF_y^{-k}\,\frac{\partial}{{\partial}y}\,\boldsymbol{I}_{n,k}\\
=\sum_j(-1)^{j+2}jF_y^{-(j+2)}F_xF_{yy}\boldsymbol{I}_{n,j}+
\sum_j(-1)^{j+1}F_y^{-(j+1)}F_x\,\frac{\partial}{{\partial}y}\,
\boldsymbol{I}_{n,j}\\
=\sum_k(-1)^k(k-2)F_y^{-k}F_xF_{yy}\boldsymbol{I}_{n,k-2}+
\sum_k(-1)^kF_y^{-k}F_x\,\frac{\partial}{{\partial}y}\,
\boldsymbol{I}_{n,k-1}.
\end{multline*}

Therefore, the derivative with respect to $x$ of the right 
side of (\ref{E:implicit}) is
\[
\sum_k(-1)^kF_y^{-k}\left[(k-1)F_{xy}\boldsymbol{I}_{n,k-1}+
\frac{\partial}{{\partial}x}\,\boldsymbol{I}_{n,k}+
(k-2)F_xF_{yy}\boldsymbol{I}_{n,k-2}+F_x\,\frac{\partial}{{\partial}y}\,
\boldsymbol{I}_{n,k-1}\right]
\]
which is
\[
\sum_k(-1)^kF_y^{-k}\boldsymbol{I}_{n+1,k}
\]
by (\ref{E:recurr}).  Therefore, (\ref{E:implicit}) holds for $n+1$
if it holds for $n$, so Theorem \ref{T:implicit} is true.

\end{document}